\documentclass[a4paper,11pt]{amsart}
\usepackage[latin1]{inputenc}
\usepackage[english]{babel}
\usepackage{latexsym}
\usepackage{amssymb}
\usepackage{amsmath}
\usepackage{amsthm}
\usepackage[all]{xy}

\newtheorem*{teos}{Theorem}
\newtheorem{teo}{Theorem}[section]
\newtheorem{lem}[teo]{Lemma}
\newtheorem{prop}[teo]{Proposition}

\theoremstyle{definition}
\newtheorem{dhef}[teo]{Definition}
\newtheorem{exe}[teo]{Example}
\newtheorem{oss}[teo]{Remark}

\newcommand{\lrg}{\longrightarrow}

\newcommand{\ze}{\mathbb{Z}}
\begin{document}
\title{Local structure of abelian covers}
\author{Donatella Iacono}
\date{2nd June 2006}
\address{\newline Dipartimento di Matematica ``Guido
Castelnuovo'',\hfill\newline Universit\`a di Roma ``La
Sapienza'',\hfill\newline P.le Aldo Moro 5, I-00185 Roma\\ Italy.}
\email{iacono@mat.uniroma1.it} \maketitle
\begin{abstract}
We study normal finite abelian covers of smooth varieties. In
particular we establish combinatorial conditions so that a normal
finite abelian cover of a smooth variety is Gorenstein or locally
complete intersection.
\end{abstract}

\section{Introduction}
Let $G$ be a finite abelian group and let $X,Y$ be algebraic
varieties. An abelian cover $X$ of $Y$, with group $G$ is a
(finite) morphism $\pi: X \lrg Y$, together with  a faithful
action of $G$ on $X$, commuting with $\pi$, such that $Y=X/G$. We
will focus our attention on the case in which $X$ is normal and
$Y$ is smooth; in this case $\pi$ is flat \cite[Sec. 3]{bib
Berry}.

\smallskip
For convenience, throughout the paper we shall  work over the
ground field $\mathbb{C}$ of the complex numbers. The results here
apply over any algebraic closed field $k$ and if $ch(k)\neq 0$ we
only need that $gcd(|G|,ch(k))=1$.

\smallskip

The theory of cyclic covers of algebraic surfaces was studied
first by A. Comessatti in \cite{bib comessatti}. Then F. Catanese
\cite{bib catanese} studied smooth abelian covers  in the case
$(\ze_2)^2$ and R. Pardini \cite{bib pardini} analyzed  the
general case. In \cite{bib nuovoManet} and \cite{bib manet} M.
Manetti investigated the property to be Gorenstein for
$(\ze_2)^n$-covers. Recently, in \cite{bib Murphy} R. Vakil used
the $(\ze_p)^n$-covers to show the existence of badly-behaved
deformation spaces.

\medskip

The first purpose of this paper is to establish combinatorial
conditions  so that a normal  $G$-cover $X$ of a smooth variety
$Y$ ($\pi:X \lrg Y$) is Gorenstein.

\smallskip

To give a more precise statement let us introduce the
combinatorial data of the cover $\pi$ at a point $y$ in $Y$. Let
$R$  be the ramification locus, i.e. the set of the points of $X$
that have non trivial stabilizer, and define the branch locus $D$
as the image of $R$ under $\pi$. By the theorem of purity of
branch locus \cite[Prop. 2]{bib zariski}, $R$ and $D$ are
divisors. We consider the divisors
with the reduced structure.

For each irreducible component T of $R$, let $H$ be the subgroup
of $G$ that stabilizes $T$; i.e. $H=\{ h \in G: hx=x \ \forall x
\in T\}$. $H$ is called the inertia group of $T$ and, by
definition, it is a finite abelian group.

Let $x$ be a smooth point in $X$ and in $T$. By Cartan's lemma
\cite[Lemma 2]{bib cartan}, the inertia group $H$ of $T$ acts
faithfully on $T_{x,X}$, leaving fixed $T_{x,T}$. Then there
exists a character $\psi: H \lrg \mathbb{C}^*$ such that the
action on $T_{x,X}/T_{x,T}$ of an element $h$ of $H$ is given by
the multiplication by $\psi(h)$. Since the action is faithful,
$\psi$ is a character generating the dual of $H$ and $H$ is a
finite subgroup of $\mathbb{C}^*$ and so it is cyclic. If $k\neq
\mathbb{C}$ we cannot use Cartan's lemma, but also in this case
$H$ is cyclic and there is an associated character $\psi$
generating the dual of $H$ \cite[Lemma 1.1 and 1.2]{bib pardini}.

Now, let $E$ be an irreducible component of $D$. Since $G$ is
abelian, all the components of $\pi^{-1}(E)$ have the same inertia
group and same associated character $\psi$. Therefore to every
component of $D$ we can associate a cyclic subgroup $H$ of $G$ and
a character $\psi $ generating the dual of $H$. For each pair
$(H,\psi)$, with $H$ a cyclic subgroup of $G$ and $\psi$ a
character generating the dual of $H$, let $D_{H,\psi}$ be the
union of all the components of $D$ that have associated  $H$ and
$\psi$.

\begin{dhef}\label{combinatorial data del punto y di Y}
Let $X$ be a normal $G$-cover of a smooth variety $Y$. Let $y$ be
a point of $Y$ lying on the components $D_{H_1,\psi_1}, \ldots,
D_{H_s,\psi_s}$ of the branch locus $D$. The set
$\{(H_i,\psi_i)\}_{i=1,\ldots,s}$ is the {\it combinatorial data
at the point $y$ in $Y$} of the cover.
\end{dhef}

\smallskip

\noindent Our first result is the following theorem (theorem
\ref{theorema A in sezion proof}) that relates the property to be
Gorenstein to the combinatorial data.
\bigskip

\begin{teos} Let X be a normal $G$-cover of a
smooth variety $Y$ and $\{(H_i,\psi_i)\}_{i=1,\ldots,s}$ its
combinatorial data at the  point $y$ in $Y$. Then the points of
$X$ over $y$ are Gorenstein points of  $X$ if and only if there
exists a character $\chi$ of $G$ such that
$$
\chi_{|H_i}=\psi_i \ \ \ \ \ \ \forall \  1 \leq i \leq s.
$$
\end{teos}
\medskip

\medskip

It is easy to see  that the existence of a character $\chi$ of $G$
is a necessary condition. Suppose that G leaves fixed the point
$x$ in $X$. Since $X$ is Gorestein the dualizing sheaf is locally
free of rank 1 \cite[Ch. 21]{bib eisenbud}. Then there exists a
character $\chi$ of $G$ such that the action of $G$ on the
dualizing sheaf is given by the multiplication by $\chi$. Let $y$
be a smooth point of $X$ near $x$. Suppose that $y$ is a smooth
point of the irreducible component $T$ of $R$ fixed by the
subgroup $H$ of $G$. As already said, we can define a character
$\psi$ of $H$, such that the faithful action of $H$ on
$T_{y,X}/T_{y,T}$ is given by the multiplication
by $\psi$. Therefore it is necessary that $\chi_{|H}=\psi$.

The previous theorem (theorem \ref{theorema A in sezion proof})
generalizes this condition and it proves that it is necessary and
sufficient condition.

\bigskip

Then we study  locally complete intersection covers. Related to
locally complete intersections there are the locally simple
covers.

\begin{dhef}\label{definiz locally simple}
Let $Y$ be a smooth variety and $\pi: X \lrg Y$ be a normal
$G$-cover. Suppose that $y \in Y$ lies on the components
$D_{H_1,\psi_1}, \ldots ,D_{H_s,\psi_s}$ of D.  The $G$-cover is
called {\it locally simple}, or $\pi$ is {\it locally simple}, if
the map $\displaystyle \oplus_i H_i \lrg G$ is injective, for each
$y \in Y$.
\end{dhef}

\noindent Every locally simple cover is a locally complete
intersection (remark \ref{oss loc siple implica sempre loc compl
inters}).
\smallskip

The second purpose of this paper is to investigate when the
converse holds.

\noindent First of all we note that the equivalence between the
properties of being locally simple and being a locally complete
intersection doesn't hold for any abelian group $G$. Actually,
example \ref{example Z_pqr cover} shows the existence of a
$G$-cover that is locally complete intersection but not locally
simple.

Hence we turn our attention on $(\ze_p)^n$-covers. For these
covers some results already exist. In \cite{bib nuovoManet}
Manetti proved that the equivalence holds for $(\ze_2)^2$-covers
\cite[Prop. 1]{bib nuovoManet}.

\noindent The following theorem (theorem \ref{teorema equiva per
(Z_p)^n}) proves that the equivalence holds for each
$(\ze_p)^n$-covers, with $p$ prime.
\bigskip

\begin{teos} Let $\pi:X \lrg Y$ be a normal flat $(\ze_p)^n$-cover with
$Y$ smooth. Then the following conditions are equivalent:
\begin{description}
   \item[i)] $\pi$ is locally simple;
   \item[ii)] X is locally complete intersection.
\end{description}
\end{teos}

This theorem is the correct version of proposition $3.25$ of
\cite{bib manet}, that states the equivalence of the previous two
properties and the property of being Gorenstein for
$(\ze_2)^n$-covers. In fact, example \ref{esempio (Z_2)^3 gor non
loc simple} shows that the property of being Gorenstein is not
equivalent to the property of being locally simple for
$(\ze_2)^n$-covers, with $n \geq 3$. However, in \cite{bib manet}
the result of proposition $3.25$ is applied only in corollary
$3.26$ where Manetti used only the true equivalence between the
property of being locally simple and being a locally complete
intersection. Therefore this minor mistake doesn't affect the rest
of the paper.  He also showed that for local simple covers the
theory of deformation is easy to understand. So locally simple
covers are interesting from the point of view of deformation
theory.

\subsection*{Notations}

  $G$ defines a finite abelian group and we use the additive
notation. In particular, we use the notation $G= \langle a \in G
:na=0\rangle$, with $n \in \mathbb{N}$, for a cyclic group $G$,
generated by $a$, of cardinality $n$. So $G \cong \ze_n$. We also
write $\langle a\rangle$ when the order is clear by the context. In general
$|G|$ is the order of a group $G$.

$G^*=\operatorname{Hom}(G,\mathbb{C}^*)$ is the group of characters  of
$G$ (the dual of $G$) and we use for it the multiplicative
notation.

$\xi \in \mathbb{C}^*$ stands for a primitive root of unity, whose
(multiplicative) order will be specified every time.

If $a,b,c \in \ze $, then $a \equiv b \ (mod \ n)$ means that $a$
is congruent to $b$ modulo $n$, $gcd(a,b)$ is the great common
divisor between $a$ and $b$ and $[a]$ is the integer part.

\bigskip

\bigskip

\noindent We collect here some results  that we will use in the
sequel.

\bigskip

\begin{teo}\label{teo matsu gor se base e fibre}
(Watanabe) If $f:X\lrg Y$ is a flat and surjective morphism of
preschemes, then $X$ is Gorenstein if
  and only if $Y$ and each fiber $f^{-1}(y)$ are Gorenstein.
\end{teo}

\begin{proof}
See \cite[Th. 1']{bib WATA INTcomple}.
\end{proof}

\begin{teo}\label{teo inters completa se base e fibre}
(Avramov) Let $f:X\lrg Y$ be a flat and surjective morphism of
locally Noetherian schemes. Then $X$ is a locally complete
intersection if and only if  the scheme $f^{-1}(y)$ is a locally
complete intersection for every point $y \in Y$ and  $Y$ is a
locally complete intersection.
\end{teo}

\begin{proof}
See \cite[Cor. 2]{bib avramov}.
\end{proof}

\begin{lem}\label{lemma hensel implica divisibile}
If $(A,m)$ is an henselian local ring with residue field $A/m =k $
algebraically closed, then the multiplicative group $A^*=A-m$ is
divisible.
\end{lem}
\noindent \begin{proof} By definition $A^*$ is divisible if and
only if for each $a \in A^*$ and for each $n \in \mathbb{N}$ there
exists $x \in A^*$ such that $a=x^n$. This condition is equivalent
to the existence of a root in $A^*$ of the monic polynomial
$P(x)=x^n-a \in A[x]$. Let $ \overline{P}(x)$ the image of $P(x)$
in $k[x]$. $ \overline{P}(x)$ has a simple root $\overline{a}$ in
$k$ and so, since $A$ is henselian,  $P(x)$ has a root $a$ that
lift $\overline{a}$ \cite[Ch. 7, Prop. 3]{bib raynaud}.

\end{proof}

\section{The key example} \label{sezio key example}

This section  is devoted to the study of a particular case of
$G$-cover $\rho:X \lrg X/G$, with $X=\mathbb{C}^s/K$ and $K$ a finite
abelian group.

\bigskip

Let $\{(H_i, \psi_i)\}_{i=1,\ldots,s}$ be a set of data, i.e.
$H_i$ are cyclic groups and  $\psi_i$ are characters generating
${H_i}^*$. If $H_i=H_j$ then we assume that $\psi_i \neq \psi_j$
in ${H_i}^*={H_j}^*$.

Let $d_i$ be the order of $H_i$, $\displaystyle
H=\bigoplus_{i=1}^s H_i$ and $K$ be a subgroup of $H$. Then we can
construct an exact sequence of abelian groups
\begin{equation}\label{o-K - H=+H_i - G - 0}
0 \lrg K \stackrel{\iota}{\lrg} H=\bigoplus_{i=1}^s H_i
\stackrel{\nu}{\lrg} G\lrg 0
\end{equation}
where $\nu$ is the sum and $G$ is the quotient. Moreover, suppose
that each $H_i$ injects into $G$.

For each element $k \in K$ we have
\begin{equation}\label{elemn k -(h_i)- sum h_i=0}
k \stackrel{\iota}{\lrg} \iota(k)=(h_1,\ldots , h_s)
\stackrel{\nu}{\lrg} \sum_{i=1}^s h_i=0.
\end{equation}
Therefore for each $k \in K-\{0\}$, $\iota(k)$ contains at least
two non-zero elements $h_i \in H_i$ and $h_j \in H_j$, for $1 \leq
i,j \leq s $.

\medskip

First of all, we note that $ H^*=\prod_{i=1}^s {H_i}^*$ and so,
using the left-exactness of the
$\operatorname{Hom}(-,\mathbb{C}^*)$ functor \cite[Prop. 2.9]{bib
Atiyah}, we have the following exact sequence
\begin{equation}\label{0-G^* - H^*-K^*}
0 \lrg G^* \stackrel{{\nu}^*}{\lrg} H^* \stackrel{{\iota}^*}{\lrg}
K^*.
\end{equation}

\bigskip

Now, we consider the action of $H$ on $\mathbb{C}^s$, with  coordinates
$(z_1,z_2, \ldots, z_s)$ given by the characters $\psi_i$ of the
  data: if $H \ni h=(h_1, \ldots , h_s)$ with $h_i \in
H_i$, then we define
$$
h(z_1,\ldots,z_s)=(\psi_1(h_1)z_1,\ldots ,\psi_s(h_s)z_s).
$$

\medskip

\noindent Each non trivial element $(0,\ldots,0,h_i,0,
\ldots,0)\in H$ acts on $\mathbb{C}^s-\{z_i=0\}$ without fixed points.
Therefore the points with all coordinates non zero have trivial
stabilizer; the hyperplanes $z_i=0$ are stabilized by $H_i\subset
H$ and they correspond to the irreducible components $R_i$ of the
ramification locus $R$.

\smallskip

Let $x \in \mathbb{C}^s$ be a point of  $R_i$. Then $H_i$ acts on the
tangent space $T_{x,\mathbb{C}^s}$ to $\mathbb{C}^s$ in $x$, leaving fixed
$T_{x,R_i}$. Since the action of $H_i$ is given by the
multiplication by the character $\psi_i$, the character associated
to the action of $H_i$ on $T_{x,\mathbb{C}^s}/T_{x,R_i}$ is $\psi_i$
itself.
\smallskip

The sub-ring of the invariant polynomials in $\mathbb{C}[z_1,\ldots,z_s]$
is the ring
$\mathbb{C}[z_1,\ldots,z_s]^H=\mathbb{C}[z_1^{d_1},\ldots,z_s^{d_s}]$
($d_i$ is the order of $H_i$), since
$h(z_i^{d_i})=(\psi_i(h_i))^{d_i} z_i^{d_i}=z_i^{d_i}$.

\noindent Set $u_i:=z_i^{d_i}$; then $\mathbb{C}[z_1,\ldots,z_s]^H=
\mathbb{C}[u_1,\ldots,u_s]$ and so $\mathbb{C}^s/H \cong \mathbb{C}^s$.

\noindent Moreover the irreducible components $D_i$ of the branch
locus $D$ are defined by the equations $u_i=0$ and they correspond
to the the components $D_{H_i,\psi_i}$ of $D$.

\bigskip

Now, we return to the exact sequence (\ref{o-K - H=+H_i - G - 0}).
Using it, we can split the action of $H$ through an action of $K$
followed by an action of $G$ and we obtain the following
commutative diagram
\begin{center} $\xymatrix{\mbox{\large\ \  $\displaystyle \mathbb{
C}^s\ar[rr]^{H\ \ } \ar[dr]_K $} &  & \mbox{\large \   $
\displaystyle \frac{\mathbb{ C}^s}{H}=\mathbb{C}^s=\displaystyle
                 \frac{X}{G} $} \\
             & \mbox{\large\ \ \ \ \  $ \displaystyle
               \frac{\mathbb{ C}^s}{K}=X.$}
               \ar[ru]_{\ \rho}^{G-cover}\  &  \\ }$
\end{center}
Since $X=\mathbb{C}^s/K$ is normal \cite[Sec. 2]{bib cartan} and
$X/G=(\mathbb{C}^s/K)/G=\mathbb{C}^s/H$ is smooth, the map
$\rho:X\lrg X/G$ is flat \cite[Sec. 3]{bib Berry}; in this way we
have obtained a $G$-cover $\rho$.

\noindent By definition, the $G$-cover $\rho$ is locally simple if
$K$ is the trivial subgroup of $H$.

\noindent Using (\ref{elemn k -(h_i)- sum h_i=0}), the action of
an element $k \in K$ is determined by the image
$\iota(k)=(h_1,\ldots h_s)$ and so each $k \in K-\{0\}$ acts non
trivially on at least two coordinates (it isn't a
pseudo-reflection).
\smallskip

As for $H$, the coordinate hyperplanes are the irreducible
components of $R$ fixed by $H_i \subset G$ (this inclusion is
guaranteed by our initial assumption $H_i$ injects in $G$).
Analogously the irreducible components of the branch locus $D$ of
$\rho$ are the $D_i$ defined by the equations $u_i=0$.

\smallskip

The coordinate ring of $X=\mathbb{C}^s/K$ is
$\mathbb{C}[z_1,\ldots,z_s]^K \subset
\mathbb{C}[z_1,\ldots,z_s]$. Since $K$ is a finite abelian group that acts
diagonally,  $\mathbb{C}[z_1,\ldots,z_s]^K$ is generated by invariant
monomials. Let $N^s_K= \{ (\alpha_1,\ldots,\alpha_s)\in
\mathbb{N}^s|\ z_1^{\alpha_1}\cdots z_s^{\alpha_s} \mbox{ is an
invariant monomial for K} \}$. Then  we can write $\displaystyle
\mathbb{C}[z_1,\ldots,z_s]^K= \bigoplus_{(\alpha_1,\ldots,\alpha_s)\in
N_K^s}z_1^{\alpha_1}\cdots z_s^{\alpha_s}\mathbb{C}$.

A monomial $z_1^{\alpha_1}\cdots z_s^{\alpha_s}$ lies in $N^s_K$
if and only if $\iota^*(\Pi_i  \psi_i^{\alpha_i})=1 \in K^*$.
Actually, each element of $\displaystyle K \leq
H=\bigoplus_{i=1}^s H_i$ is of the form $k=(k_1,\ldots,k_s)$ and
so $k(z_1^{\alpha_1}\cdots z_s^{\alpha_s})=z_1^{\alpha_1}\cdots
z_s^{\alpha_s}$ ($\forall \ k$) if and only if $
\psi_1(k_1)^{\alpha_1} z_1^{\alpha_1} \cdots
\psi_s(k_s)^{\alpha_s} z_s^{\alpha_s} = z_1^{\alpha_1} \cdots
z_s^{\alpha_s}$ that is $\Pi_i  \psi_i^{\alpha_i} =1 \in K^*$.

By the exact sequence (\ref{0-G^* - H^*-K^*}), this is equivalent
to the existence of a character $\chi \in G^*$ such that
$\nu^*(\chi)=\Pi_i  \psi_i^{\alpha_i}$. We notice that the natural
numbers $\alpha_i$ are such that $\chi_{|H_i}= \psi_i^{\alpha_i}$,
for each $i$.

Therefore $\displaystyle \mathbb{C}[z_1,\ldots,z_s]^K=
\bigoplus_{(\alpha_1,\ldots,\alpha_s)\in N_K^s}
z_1^{\alpha_1}\cdots z_s^{\alpha_s}\mathbb{C} = \bigoplus_{\chi \in
G^*}\mathcal{A}_\chi$ where the $\mathcal{A}_\chi$ are defined as
follows.

\noindent Since the groups $H_i^*$ have order $d_i$,
$\mathcal{A}_1=\mathbb{C}[z_1^{d_1},\ldots,z_s^{d_s}]=
\mathbb{C}[u_1,\ldots,u_s]$; moreover  $\mathcal{A}_\chi$ are free
$\mathcal{A}_1$-module of rank 1. In fact, if the character
$\chi=\prod_i \psi_i^{\alpha_i} \in G^*$  is associated to the
invariant monomial $z_1^{\alpha_1} \cdots z_s^{\alpha_s}$, then
the other invariant monomials $z_1^{\beta_1} \cdots z_s^{\beta_s}$
with the same associated character $\chi=\prod_i
\psi_i^{\alpha_i}$ are such that $\beta_i \equiv \alpha_i \ (mod \
d_i)$, for each $0 \leq i \leq s$. Hence, if $z_1^{\alpha_1}
\cdots z_s^{\alpha_s}$ is the unique invariant monomials
associated to $\chi$, with $0 \leq \alpha_i \leq d_i-1$, we write
$w_\chi=z_1^{\alpha_1} \cdots z_s^{\alpha_s}$ and so $w_\chi$ is a
generator of $\mathcal{A}_\chi$ as free $\mathcal{A}_1$-module,
i.e. $\mathcal{A}_\chi=\langle w_\chi=z_1^{\alpha_1} \cdots
z_s^{\alpha_s}\rangle$.

\medskip

Now we analyze the multiplicative structure of $\mathbb{C}[z_1,
\ldots,z_s]^K$.\\
Given $w_\chi=z_1^{\alpha_1} \cdots z_s^{\alpha_s}$ and
$w_{\chi'}=z_1^{\alpha_1'} \cdots z_s^{\alpha_s'}$, we define

\begin{equation}\label{epsilo MIO uguale epsilon building data}
\varepsilon^i_{\chi,\chi'} :=[\frac{\alpha_i+\alpha_i'}{d_i}].
\end{equation}
As already observed, all the $\alpha_i$ and $\alpha_i'$ satisfy
$\chi_{|H_i}= \psi^{\alpha_i}$ and  $\chi'_{|H_i}=
\psi^{\alpha_i'}$.

\noindent Fix the attention on the fiber of $\rho$ over $0 \in
\mathbb{C}^s/H$: $0$ lies on the components $D_{H_i,\psi_i}$ for each
$i=1,\ldots,s$ and so $u_i=0$. Then $w_\chi \cdot w_{\chi'}=0$ if
there exists $i$ such that $\alpha_i+\alpha_i' \geq d_i$,
otherwise $w_\chi \cdot w_{\chi'}=w_{\chi \cdot \chi'}$.

\noindent Hence the product in $\mathbb{C}[z_1,\ldots,z_s]^K$ is given by
\begin{equation}\label{equazio MIA BUILD DATA}
w_\chi \cdot w_{\chi'}=z_1^{\alpha_1+\alpha_1'}\cdots
z_s^{\alpha_s+\alpha_s'}=w_{\chi \cdot \chi'}\prod_{i=1}^s u_i^{
\varepsilon^i_{\chi,\chi'}}.
\end{equation}

\subsection{When is the key example Gorenstein?}

Now we investigate when $X$ of the $G$-cover $\rho:\mathbb{C}^s/K=X\lrg
X/G=\mathbb{C}^s/H=\mathbb{C}^s$ is Gorenstein, establishing combinatorial
conditions on the data $\{(H_i, \psi_i)\}_{i=1,\ldots,s}$.

\medskip

Since $X=\mathbb{C}^s/K$, we can apply an useful result due to
Watanabe \cite[Sec. 2, Th. 1]{bib watanabe}.
\begin{teo}\label{teo watanabe}(Watanabe)
Let $R=k[x_1,\ldots,x_n]$ be a polynomial ring over a field $k$
and $G$ be a finite subgroup of $GL(n,k)$ with $gcd(|G|,ch(k))=1$
if $ch(k)\neq 0$. We also assume that $G$ contains no pseudo
reflections. Then $R^G$ is Gorenstein if and only if $G \leq
SL(n,k)$.
\end{teo}

In our example the action of $K$ just depends on the action of $H$
and so it is strictly related on the data $\{(H_i, \psi_i)\}_{
i=1,\ldots,s}$. Therefore we can rewrite the previous theorem
\ref{teo watanabe} in term of combinatorial conditions on the
characters $\psi_i$.

More precisely, we have associated  to each element $k \in K$ a
diagonal matrix of the action on $\mathbb{C}^s$, whose entries are
the characters $\psi_i$ valuated on $k$. Therefore we need to know
when a product of characters $\psi_i$ is $1$ for each element $k
\in K$. Using the exact sequence (\ref{0-G^* - H^*-K^*}), the
previous condition is related to the existence of a character of
$G$ that lifts the characters $\psi_i$. In conclusion we have
proved the following proposition that is a particular case of the
first  theorem.

\begin{prop}\label{pro caso particolare per gruppi di teo A}
Let $\{(H_i,\psi_i) \}_{i=1,\ldots ,s}$ be a set of data and K a
subgroup of $\displaystyle H=\bigoplus_{i=1}^s H_i$. Let $G$ be
the quotient and suppose that each $H_i$ injects in $G$. $\mathbb{C}^s/K$
is Gorenstein if and only if there exists a character $\chi$ of
$G$ such that
\begin{equation}\label{equa chi|=psi_i in proposizio del teo A}
\chi_{|H_i}=\psi_i \ \ \ \ \ \ \forall \ 1 \leq i \leq s.
\end{equation}
\end{prop}

\subsection{Examples}\label{sottosezio esempi}

Now, we explicitly analyze some $G$-covers. Fixing a finite
abelian group $G$, we can consider a set of combinatorial data
$\{(H_i,\psi_i) \}_{i=1,\ldots ,s}$, with $H_i$ cyclic subgroups
of $G$ (so each one injects in $G$) and $\psi_i$ a character
generating $H_i^*$. If $H_i=H_j$ then we assume that $\psi_i \neq
\psi_j$ in ${H_i}^*={H_j}^*$. Suppose that $\nu:
H=\bigoplus_{i=1}^s H_i \lrg G$ is surjective. Also in this case
we can construct an exact sequence as (\ref{o-K - H=+H_i - G - 0})
with  $K$ the kernel of $\nu$.

\bigskip

\begin{exe} \label{esempio (Z_2)^3 gor non loc simple}
Let $G \cong (\ze_2)^3$ with standard generators $e_i$ and let
$H_i=\langle e_i\rangle$, for $i=1,2,3$, and $H_4=\langle e_1+e_2+e_3\rangle$ with
associated character $\psi_1,\ \psi_2,\ \psi_3$ and $\varphi$,
respectively ($\psi_i(e_i)=\xi$, $\varphi(e_1+e_2+e_3)=\xi$, with
$\xi \in
\mathbb{C}^*$ and $\xi^2=1$).\\
In this example, $H=(\ze_2)^4$  and for $k \in K-\{0\}$ we have
$\iota(k)=(e_1,e_2,e_3 ,e_1+e_2+e_3)\in H$. The diagram associated
to the action of $H$ on $\mathbb{C}^4$ is the following
\begin{center}
$\xymatrix{\mbox{\large\ \  $\displaystyle \mathbb{
C}^4\ar[rr]^{H} \ar[dr]_K $} &  &
  \mbox{\large\ \  $\displaystyle \frac{\mathbb{
C}^4}{H}=\frac{X}{(\ze_2)^3} $}  \\
             & \mbox{\large  \ \ \ \ \ $\displaystyle  \frac{\mathbb{ C}^4}
             {\pm Id} $}=X \ar[ru]_{\ \ (\ze_2)^3-cover} &  \\
             }$
\end{center}
with actions given by:\\
- $h \in H$ is of the form $h=t_1\, e_1+ {t_2 }\, {e_2}+{t_3 }\,
{e_3}+t_4\,( e_1+e_2+e_3)$, with $t_i=0,1$; so
$h(z_1,z_2,z_3,z_4)=((-1)^{t_ 1}z_1 ,(-1)^{t_2 }z_2,(-1)^{t_3 }z_3
,(-1)^{t_4 }z_4 )$.\\
- $k \in K-\{0\}$ acts as: $k(z_1,z_2,z_3,z_4)=
(-z_1,-z_2,-z_3,-z_4)$.\\
- $g \in G$ is of the form $g=s_1 \,e_1+ s_2 \, e_2+ s_3\, e_3$,
with $s_i=0,1$; so $g(w_1,w_2,w_3,w_4)=((-1)^{s_ 1}w_1 ,(-1)^{s_2
}w_2,(-1)^{s_3 }w_3 ,w_4 )$.

\smallskip

The $(\ze_2)^3$-cover $\rho: X \lrg X/(\ze_2)^3$ is not locally
simple $(K \neq \{0\})$ but,  since $K\leq SL(4,\mathbb{ C})$, $X$
is Gorenstein by theorem \ref{teo watanabe}.

\begin{oss}\label{osserv(Z_2)^3-cover no loc sem no inters compl}
In this example X cannot be a local complete intersection.\\
Let us recall a theorem due to Schlessinger \cite[Sec. 10]{bib
Stevens}  or \cite[Sec. 3]{bib schlessin}.
\begin{teo}\label{teorema schlessinger} (Schlessinger)
Quotient singularities which are non singular in codimension two
are rigid, that is every infinitesimal deformation is trivial.
\end{teo}

\noindent We can apply this theorem to the cover $\mathbb{C}^4
\lrg \mathbb{C}^4/K$. The quotient $\mathbb{C}^4/K$ has a
singularity only at the origin ($0$ is the only fixed point) and
its codimension is 4. Therefore every infinitesimal deformation of
$X=\mathbb{C}^4/K$ is trivial. On the other hand, singular
complete intersections are not infinitesimally rigid (cf.
\cite{bib Artin}). Therefore $X$ cannot be a local complete
intersection.
\end{oss}
\end{exe}

\begin{exe}\label{example Z_pqr cover}
Let $G \cong \ze_{pqr}$, with $p< q <r$ prime numbers. Let $\xi
\in \mathbb{C}^*$ be a primitive root of $1$ of order $pqr$, and
$G$ generated by $a$. Let $H_1 \cong \ze_{pr}=\langle q \,
a\rangle $ and $H_2 \cong \ze_{pq}=\langle r \, a\rangle$,
$\psi_1(qa)=\xi^{q \alpha}$, with $gcd(\alpha,pr)=1$, and
$\psi_2(ra)=\xi^{r\beta}$, with $gcd(\beta,pq)=1$. In this
example, $H=H_1 \oplus H_2$ and $K\cong \ze_p$; so we have the
following commutative diagram
\begin{center}
$\xymatrix{\mbox{\large\ \  $\displaystyle \mathbb{ C}^2\ar[rr]^{\
H} \ar[dr]_K $}  &  &
  \mbox{\large   \  $\displaystyle \frac{\mathbb{ C}^2}{H} $}  \\
             & \mbox{\large \ \ $\displaystyle  \frac{\mathbb{ C}^2}{K}. $}
  \ar[ru]_{\ze_{pqr}-cover} &  \\
             }$
\end{center}
Applying proposition \ref{pro caso particolare per gruppi di teo
A}, $\mathbb{C}^s/K$ is Gorenstein  if and only if ${\psi_1}_{|K}
={\psi_2}_{|K} $ or equivalently $\alpha \equiv \beta\ (mod\, p)$.

Moreover, if $\psi_1$ and $\psi_2$ satisfy this condition then
$\mathbb{C}^2/K$ is a local complete intersection
$\ze_{pqr}$-cover that isn't locally simple ($K\cong \ze_p$). In
fact, if $\alpha \equiv \beta\ (mod \ p)$ the action of an element
$k \in K$ on $\mathbb{C}^2$, with coordinates $(z_1,z_2)$, is of
the form $k(z_1,z_2)=(z_1\eta,z_2\eta^{-1})$, where $\eta \in
\mathbb{C}^*$ is a $p$-root of 1. Therefore $\mathbb{C}^2/K$ has a
rational double point at $0$ (of type $A_{p-1}$) and it is a local
complete intersection.
\end{exe}

\begin{exe}\label{esempio Z_(p^n) ciclico}
Let  $G\cong \ze_{p^n}$, with $p$ prime. All the subgroups of $G$
are cyclic $p$-groups.

\noindent Let $\{(H_i ,\psi_i)\}_{i=1,\ldots , s}$ be a set of
combinatorial data. Then $H=\oplus_iH_i$ is the sum of $p$-groups
with at least one isomorphic to $G$ (to have a surjection).
Relabelling if necessary the combinatorial data, we can suppose
that $H=H_1\oplus \cdots \oplus H_s$, with $H_s\cong G$. Therefore
$\psi_s$ is a character of $G$ and so, according to proposition
\ref{pro caso particolare per gruppi di teo A}, the condition for
$X=\mathbb{C}^s/K$ to be Gorenstein becomes
${\psi_s}_{|H_i}=\psi_i$, for each $1 \leq i \leq s-1$. Moreover,
the group $G$ is cyclic and so there exists only one cyclic
subgroup of $\ze_{p^n}$ of order $p^i$, for each $i=1,\ldots, n$.
For this reason the cyclic groups $H_i$ must have different orders
and so $s \leq n$. Then we can reorganize them so that $H_1 < H_2
< \cdots < H_s =G$ and  the following restrictions must match,
i.e: $\psi_s|_{H_j}=\psi_j=\psi_{j+1}|_{H_j}$, for each $1 \leq j
\leq s-1$.
\end{exe}

\begin{exe}\label{esempio (Z_p)^n-cover con remark n_i mag 3}
Let $G\cong  {(\ze_p)}^n$ ($p$ prime) and  $\{(H_i
,\psi_i)\}_{i=1,\ldots , s}$ be a set of combinatorial data. In
this case the only admissible cyclic subgroups of $G$ have order
$p$ ($H_i\cong \ze_p$); therefore the exact sequence (\ref{o-K -
H=+H_i - G - 0}) can be written as
\begin{equation}
0 \lrg K \lrg H=\bigoplus_{i=1}^s \ze_p \lrg G \lrg 0
\end{equation}
with $s \geq n$.

\begin{oss}\label{oss se gorens n_i div 1}
If $s=n$ then $H=G$ and we have a locally simple cover.

If $ s >n$ and the ${(\ze_p)}^n$-cover $\rho: \mathbb{C}^s/K= X \lrg
X/{(\ze_p)}^n$ is Gorenstein, then $H_i \neq H_j$ for each $1 \leq
i,j \leq s$ and $i \neq j$. In fact, assume that $H_i =H_j$; then
by definition $\psi_i \neq \psi_j$. If $X$ is Gorenstein then
there exists a character $\chi$ of $G$ such that
$\psi_i=\chi_{|H_i}= \chi_{|H_j}=\psi_j$ that is a contradiction.
In particular $H_i \cap H_j= \{0\}$.
\smallskip

Therefore, looking at the sequence (\ref{elemn k -(h_i)- sum
h_i=0}), for each element $k \in K-\{0\}$ we have $\iota
(k)=(h_1,\ldots , h_s)$ with at least three elements $h_j$
different from zero, with $1 \leq j \leq s$.

\noindent This implies that each element $k \in K-\{0\}$ acts  non
trivially at least on three coordinates of the points in $\mathbb{C}^s$.

This fact will  be used in the proof of the main theorem (theorem
\ref{teorema equiva per (Z_p)^n}).
\end{oss}
\end{exe}

\section{Structure of abelian cover}\label{sezione preliminari}

In \cite{bib pardini}, Pardini completely described the structure
of normal abelian covers of smooth complete algebraic varieties
(theorem \ref{teo Pardini}); she also  gave conditions so that a
normal $G$-cover of a smooth variety is smooth (theorem \ref{teo
pardini liscio}). In this section we introduce some notations to
recall these theorems. A detailed description can be found in
\cite{bib pardini}.

\medskip

Let $\pi:X \lrg Y$ be an abelian $G$-cover with $X$ normal and $Y$
smooth; then $\pi$ is flat and $\pi_*\mathcal{O}_X$ is locally
free \cite{bib Berry}. The action of $G$ on $X$ induces the
splitting $\displaystyle \pi_*\mathcal{O}_X=\bigoplus_{\chi \in
G^*}L_\chi^{-1}$, with $L_\chi^{-1}$ line bundle, on which $G$
acts via the character $\chi$, and $L_1$  isomorphic to
$\mathcal{O}_Y$.

\smallskip

Let $D_{H,\psi}$ be the union of all the components of $D$ that
have associated the same subgroup $H$ and  character $\psi$.
$L_\chi$ and  $D_{H,\psi}$ are called the {\it building data  of
the cover}.

\smallskip

Suppose $\chi, \chi' \in G^*$. Then $\chi_{|H}$ and ${\chi'}_{|H}$
belong to $H^*$ and so there exist $i_ \chi$ and $i_{\chi'}$ in
$\{ 0,1, \ldots , |H|-1\}$, such that $\chi_{|H}=\psi^{i_ \chi}$
and ${\chi'}_{|H}=\psi^{i_ {\chi'}}$. Finally, let
$\varepsilon_{\chi,\chi'}^{H,\psi}= [\frac{{i_ \chi}+i_
{\chi'}}{|H|}]$. \\
Using the above notations, we can state the following theorem due
to Pardini.
\begin{teo}\label{teo Pardini}(Pardini)
Let $G$ be an abelian group. Let Y be a smooth variety, $X$ a
normal one and let $\pi: X \lrg Y$ be an abelain cover with group
$G$. Then the following set of linear equivalences is satisfied by
the building data of the cover:
\begin{equation}\label{eqn pardini building data}
L_\chi + L_{\chi'} = \sum_{H,\psi}
\varepsilon_{\chi,\chi'}^{H,\psi} D_{H,\psi}
\end{equation}
Conversely, to any set of data $L_\chi,\ D_{H,\psi}$ satisfying
(\ref{eqn pardini building data}) we can associate an abelian
cover  $\pi: X \lrg Y$ in a natural way. Whenever the cover so
constructed is normal, $L_\chi$ and $\ D_{H,\psi}$ are its
building data.\\
Moreover, if $Y$ is complete, then the building
data determine the cover $\pi :X \lrg Y$ up to isomorphism of
Galois cover.
\end{teo}
\begin{proof}
See \cite[th. 2.1]{bib pardini}.
\end{proof}

\begin{oss}
About the existence of the cover, let $L_\chi$ and $D_{H,\psi}$
($L_1=\mathcal{O}_Y$) satisfy relations (\ref{eqn pardini building data}).
Let $\sigma_{H,\psi} \in \mathcal{O}(D_{H,\psi})$ be a section defining
$D_{H,\psi}$ and $\displaystyle A=\bigoplus_{\chi \in
G^*}L_\chi^{-1}$. Using (\ref{eqn pardini building data}), the
formula
\begin{equation}\label{equaz mu_x,x' su algebra A}
\mu_{\chi,\eta}=\prod_{H,\psi}
{\sigma_{h,\psi}}^{\varepsilon_{\chi,\eta}^{H,\psi}}
\end{equation}
defines an associative multiplication on $A$. $G$ acts on each
$L_\chi^{-1}$ by multiplication by $\chi$ and so  we can extend
this action to an  action of $G$ over $A$, compatible with the
multiplicative structure. Therefore $X=\operatorname{Spec} A$ is a
$G$-cover of $Y$.

The equation (\ref{equazio MIA BUILD DATA}) of section \ref{sezio
key example} is the explicit version of (\ref{equaz mu_x,x' su
algebra A}) for the key example.
\end{oss}

\begin{prop}\label{prop su etale buildin determina cover}
Let $\pi:X \lrg Y$ and $\pi':X \lrg Y$ be $G$-covers with the same
building data. For each point $y \in Y$ there exists an \'etale
neighborhood $U_y$ over which the covers $\pi$ and $\pi'$ are
isomorphic.
\end{prop}
\noindent \begin{proof} Assume that $\pi:X \lrg Y$ and $\pi':X
\lrg Y$ are two $G$-covers with the same building data. Let
$\displaystyle A=\bigoplus_{\chi \in G^*}L_\chi^{-1}$ and denote
with $\mu_{\chi,\eta}$ and $\mu'_{\chi,\eta}$ the two algebra
structures on $A$ corresponding to $\pi$ and $\pi'$. Since the
relations (\ref{eqn pardini building data}) must be satisfied, by
(\ref{equaz mu_x,x' su algebra A}) we can conclude that
$\mu_{\chi,\eta}$ and $\mu'_{\chi,\eta}$ have the same divisors,
for each $\chi,\eta \in G^*$; therefore, for each $\chi,\eta \in
G^*$, there exists a section $c_{\chi,\eta} \in
\Gamma(Y,\mathcal{O}_{Y}^*)$ such that
\begin{equation}\label{equa mu = c_xx' mu'}
\mu_{\chi,\eta}=c_{\chi,\eta}\mu'_{\chi,\eta}.
\end{equation}
Since $\mu_{\chi,\eta}$ and $\mu'_{\chi,\eta}$ are associative
multiplications, the sections $c_{\chi,\eta}$ satisfy the
following identity
\begin{equation}\label{equa associativi c_xx'}
c_{\eta,\tau}c_{\eta\tau,\chi}=c_{\chi,\eta}c_{\chi\eta,\tau} \ \
\ \ \ \forall \chi,\eta,\tau \in G^*.
\end{equation}
Let $\varphi_\chi$ be automorphisms of the invertible sheaves
$L_\chi$. If $w_\chi$ generates $L_\chi$ then there exists $a_\chi
\in \Gamma(Y,\mathcal{O}_{Y}^*)$ such that $\varphi_\chi(w_\chi)=a_\chi
w_\chi$.

By definition a $G$-isomorphism $\varphi:(A,\mu) \lrg (A,\mu')$ of
the $G$-algebras $A$ with multiplication $\mu$ and $A$ with
$\mu'$, is such that $\varphi(w_\chi \cdot w_\eta) =
\varphi(w_\chi) \cdot \varphi(w_\eta)$. This condition is
satisfied if and only if $a_{\chi\eta}\mu_{\chi,\eta} w_{\chi
\eta}=a_\chi a_\eta \mu'_{\chi,\eta}w_{\chi \eta}$. Therefore,
using the identity (\ref{equa mu = c_xx' mu'}), to produce a
$G$-isomorphism $\varphi$ we have
  to show the existence of elements $a_\chi \in
\Gamma(Y,\mathcal{O}_{Y}^*)$, for each $\chi \in G^*$, such that
\begin{equation}\label{equa c_xx'=a_x a_x'/a_xx'}
c_{\chi,\eta}=\frac{a_\chi a_\eta}{a_{\chi\eta}} \ \ \ \ \ \
\forall \ \chi, \eta\in G^*.
\end{equation}

Fix the attention over a point $y \in Y$. To prove the existence
of the elements $a_\chi$ we use the cohomology of the group $G^*$
with coefficients in $B=\mathcal{O}_{y,Y}^*$ (with multiplicative notation)
considered as trivial $G^*$-module.

The elements $c_{\chi,\eta}$ can be considered as elements of
$C^2(G^*,\mathcal{O}_{y,Y}^*)$ by $c(\chi,\eta)=c_{\chi,\eta} \in
\mathcal{O}_{y,Y}^*$. The associative relations (\ref{equa
associativi c_xx'}) for $c_{\chi,\eta}$ is exactly the 2-cocycle
condition for $c({\chi,\eta})$  (see \cite[Ch. VII, sec. 3]{bib
Serre}). Therefore $c_{\chi,\eta}\in
H^2(G^*,\mathcal{O}_{y,Y}^*)$.

In general $H^2(G^*,\mathcal{O}_{y,Y}^*)\neq 0$. Let $\widetilde{\mathcal{O}_{y,Y}}$
be the henselianization of $\mathcal{O}_{y,Y}$. By lemma \ref{lemma hensel
implica divisibile}, $\widetilde{\mathcal{O}_{y,Y}}^*$ is divisible and so
$H^q(G^*,\widetilde{\mathcal{O}_{y,Y}}^*)=0$ for each $q \geq 1$. This
implies that there exists $h \in C^1(G^*,\widetilde{\mathcal{O}_{y,Y}}^*)$
such that $c= d h$, or better
$c_{\chi,\eta}=\frac{h(\chi)h(\eta)}{h(\chi \eta)}$ for each
$\chi,\eta \in G^*$.

The henselianization $\widetilde{\mathcal{O}_{y,Y}}$  of
$\mathcal{O}_{y,Y}$, is an inductive limit of local rings
$(B_i,m_i)$ with $B_i$ \'etale over $\mathcal{O}_{y,Y}$ (see
\cite[Th. 1, pag. 87]{bib raynaud}). Then there exists a local
ring $(B,m)$ with $B$ \'etale over $\mathcal{O}_{y,Y}$, such that
$h \in C^1(G^*,B^*)$. Let $a_\chi=h(\chi)$ for each $\chi \in
G^*$. Then the identity (\ref{equa c_xx'=a_x a_x'/a_xx'}) holds
and so in the \'etale neighborhood $U=\operatorname{Spec} B$ of
$y$ the two cover $\pi$ and $\pi'$ are isomorphic.
\end{proof}

\begin{oss}\label{oss fibra dipende solo combinator data}
By theorem \ref{teo Pardini} and proposition \ref{prop su etale
buildin determina cover},  the fiber of a $G$-cover $\pi:X \lrg Y$
over a point $y \in Y$ just depends on the combinatorial data at
$y$ of $\pi$.
\end{oss}

\bigskip

Now suppose that $L_\chi$ are invertible sheaves and $D_{H,\psi}$
are divisors, on a variety $Y$, satisfying  conditions (\ref{eqn
pardini building data}).  Let $\mathfrak{M}/\mathfrak{M}^2$ be the
cotangent space. The following theorem, due to Pardini, determines
the conditions on the building data so that an  abelian cover of a
smooth variety is smooth.

\begin{teo}(Pardini)\label{teo pardini liscio}
Assume that $Y$ ia a smooth variety and that $\pi :X \lrg Y$ is
the $G$-cover of $X$ with building data $L_\chi$'s and
$D_{H,\psi}$'s given by the previous theorem \ref{teo Pardini}.
Then $X$ is non singular over a point $y \in Y$ if and only if one
of the following conditions holds:

\begin{description}
\item[i)] y is not a branch point of $\pi$;
\item[ii)] y belongs to one component $\Delta$ of D and y is a smooth
point of $\Delta$;
\item[iii)] y lies on components $D_{H_1,\psi_1}, \ldots ,D_{H_r,\psi_r}$
of D and:
\begin{description}
   \item[a)] the map $H_1 \oplus \cdots \oplus H_r \lrg G$ is an
   injection,
   \item[b)] let $b_i$ be a local equation for $D_{H_i,\psi_i}$
   around y, $i=1,\ldots,r$. Then the subspaces of
   $\mathfrak{M}/\mathfrak{M}^2$ generated by $db_1, \ldots, db_r$ has
   dimension r.
   \end{description}
\end{description}
\end{teo}
\begin{proof}
See \cite[Prop. 3.1]{bib pardini}.
\end{proof}

\begin{oss}\label{oss loc siple implica sempre loc compl inters}
Pardini proved that  a locally simple cover is locally complete
intersection \cite[Prop. 2.1]{bib pardini}. Theorem \ref{teorema
equiva per (Z_p)^n} proves that the converse is true for
$(\ze_p)^n$-covers, for each prime $p$ and natural number $n$ and
example \ref{example Z_pqr cover} proves that it isn't true in
general.
\end{oss}

\section{Proof of the first theorem (general case)}\label{sezio proof
of theorem A}

In this section we prove the theorem that relates the property to
be Gorenstein to the combinatorial data of a cover. The main idea
is to reduce the general case of a $G$-cover $\pi: X \lrg Y$ to a
$G$-cover of the form $\rho :\mathbb{C}^s/K \lrg
(\mathbb{C}^s/K)/G$ and then apply the proposition \ref{pro caso
particolare per gruppi di teo A} of the previous section.

\bigskip

Let $\pi:X \lrg Y$ be a  $G$-cover of algebraic varieties  with
$X$ normal and $Y$ smooth. Since $G$ is a finite group, the fibers
are $0$-dimensional and their cardinalities divide the order of
$G$. The property to be Gorenstein is local and so we restrict our
attention on $X$ over a fixed point $y$ in $Y$.

By theorem \ref{teo matsu gor se base e fibre}, $X$ is Gorenstein
over $y$ if and only if the fiber of the $\pi$-cover over $y$ is
Gorenstein; therefore, from now on, we will restrict our attention
just on the fibers.

Suppose that $y \in Y$ lies  on the components
$D_{H_1,\psi_1},\ldots ,D_{H_s,\psi_s}$ of the branch locus $D$.
Therefore the set $\{(H_i,\psi_i)\}_{i=1,\ldots,s}$ is the
combinatorial data at the point $y \in Y$ of $\pi$ and, by remark
\ref{oss fibra dipende solo combinator data}, it determines the
fiber of $\pi$ over $y$.

Let  $H$ be the direct sum of the cyclic groups $H_i$, i.e.
$\displaystyle H=\bigoplus_{i=1}^s H_i$. Then we have a map $\nu$
from $H$ to $G$ ($\nu$ is the sum).

It is not restrictive to assume that $\nu$ is a surjection, since
it is possible to factorize the cover $\pi$ near $y$ as the
composition of a totally ramified cover (i.e. inertia groups
generate $G$) followed by an \'etale map in the following way. Let
$M$ be the subgroup of $G$ generated by the inertia groups and $T$
the quotient $G/M$. Then we have the commutative diagram
\begin{center}
$\xymatrix{\mbox{\large\ \  $\displaystyle X \ar[rr]^{G\ \ }
\ar[dr]_M $} &  &
\mbox{\large \  $ \displaystyle \frac{X}{G} = Y = \frac{Z}{T}$} \\
             & \mbox{\large \ \ \  $ \displaystyle \frac{X}{M}=Z.$}
\ar[ru]_T &  \\ }$
\end{center}
where $X \lrg Z$ is a totally ramified cover and $Z\lrg
\displaystyle \frac{Z}{T}$ is \'etale \cite[pag. 487]{bib
catanese}.

\smallskip

If $K$ is the kernel of $\nu$, as in section \ref{sezio key
example}, we can consider the exact sequence of abelian groups
$$
0 \lrg K \stackrel{\iota}{\lrg} H=\bigoplus_{i=1}^s H_i
\stackrel{\nu}{\lrg} G\lrg 0.
$$
and the commutative diagram
\begin{center}
$\xymatrix{\mbox{\large\ \  $\displaystyle \mathbb{
C}^s\ar[rr]^{H\ \ } \ar[dr]_K $} &  &
\mbox{\large \  $ \displaystyle \frac{\mathbb{C}^s}{H}=\mathbb{C}^s $} \\
             & \mbox{\large \ \  $ \displaystyle \frac{\mathbb{
C}^s}{K}.$} \ar[ru]^\rho_{G-cover} &  \\ }$
\end{center}
Since $H_i$ are cyclic subgroups of $G$, they  inject in $G$. In
this way we have obtained another $G$-cover $\rho$.

\medskip

Now we prove a fundamental lemma for the proof of the theorem.

\begin{lem}\label{lemm esiste combina cover stesse fibre}
The fiber of the $\pi$-cover over  $y \in Y$ is isomorphic to the
fiber of the $\rho$-cover over $0 \in \mathbb{C}^s/H$.
\end{lem}
\noindent \begin{proof} Let us consider the $G$-cover $\pi:X\lrg Y$.
By hypothesis,  $y \in Y$ lies on the components
$D_{H_1,\psi_1},\ldots ,D_{H_s,\psi_s}$ of the branch locus $D$.
Therefore $\{(H_i,\psi_i)\}_ {i=1,\ldots,s}$  is the set of the
combinatorial data at the point $y \in Y$ of $\pi$; moreover, by
theorem \ref{teo Pardini}, the building data of this cover satisfy
the relation (\ref{eqn pardini building data}) with
$\varepsilon_{\chi,\chi'}^{H,\psi}=
\varepsilon_{\chi,\chi'}^{H_j,\psi_j}= [\frac{{i_ \chi}+i_
{\chi'}}{|H_j|}]$.

Now, let us examine the $G$-cover $\rho: \mathbb{C}^s/K \lrg \mathbb{C}^s/H $
constructed above. We analyzed this type of cover in section
\ref{sezio key example}. In this case $0 \in \mathbb{C}^s/H$ lies on the
components $D_j=D_{H_j,\psi_j}$, for $j=1,\ldots,s$ and so the
combinatorial data $\{(H_j,\psi_j)\}$ of $\rho$ at $0 \in \mathbb{C}^s/H$
is the same of the $\pi$-cover in $y$. Moreover, the building data
of the $\rho$-cover satisfy equation (\ref{eqn pardini building
data}) with coefficient $ \varepsilon^j_ {\chi,\chi'}$ and by
definition $\varepsilon_{\chi,\chi'}^{H_j,\psi_j}=\varepsilon^j_
{\chi,\chi'}$.

Applying remark \ref{oss fibra dipende solo combinator data} we
can conclude that the fiber of the $\rho$-cover in $0$ is
isomorphic to the fiber of the $\pi$-cover in $y$.
\end{proof}

\begin{dhef}
Let $\pi: X \lrg Y$ be a $G$-cover with $X$  normal and $Y$ smooth
variety. The $\rho$ cover defined above is called the { \it
combinatorial cover associated to $\pi$ in $y$}.
\end{dhef}

Now we can prove the theorem.

\begin{teo}\label{theorema A in sezion proof}
Let X be a normal $G$-cover of a smooth variety $Y$ and
$\{(H_i,\psi_i)\}_{i=1,\ldots,s}$ its combinatorial data at the
point $y$. Then the points of $X$ over $y$ are Gorenstein points
of  $X$ if and only if there exists a character $\chi$ of $G$ such
that
\begin{equation}\label{cond su caratteri in teoA in sez proof}
\chi_{|H_i}=\psi_i \ \ \ \ \ \ \forall \  1 \leq i \leq s.
\end{equation}
\end{teo}
\noindent \begin{proof} By theorem \ref{teo matsu gor se base e
fibre}, $X$ is Gorenstein over $y$ if and only if the fiber
$\pi^{-1}(y)$ is Gorenstein. Let $\rho$ be the combinatorial cover
associated to $\pi$ in $y$. Applying the previous lemma \ref{lemm
esiste combina cover stesse fibre}, the fiber $\pi^{-1}(y)$ is
Gorenstein if and only if the fiber $\rho^{-1}(0)$ is Gorenstein.
Applying again theorem \ref{teo matsu gor se base e fibre},
$\rho^{-1}(0)$ is Gorenstein if and only if $\mathbb{C}^s/K$ is Gorenstein
over 0. Now we can apply proposition \ref{pro caso particolare per
gruppi di teo A} and conclude the proof.
\end{proof}

\section{Proof of the main theorem}\label{sezione proof of theorem B}

In this section  we  prove the equivalence between the properties
of being locally simple and being a locally complete intersection
for $(\ze_p)^n$-covers, for each prime $p$ and natural numbers
$n$.

\begin{teo}\label{teorema equiva per (Z_p)^n}
Let $\pi:X \lrg Y$ be a normal flat $(\ze_p)^n$-cover with $Y$
smooth. Then the following conditions are equivalent:
\begin{description}
   \item[i)] $\pi$ is locally simple;
   \item[ii)] X is locally complete intersection.
\end{description}
\end{teo}
\begin{proof} i) $\Rightarrow$ ii). It follows from
definition \ref{definiz locally simple} of locally simple and from
remark
\ref{oss loc siple implica sempre loc compl inters}.\\
ii) $\Rightarrow$ i). Suppose that the cover is not locally
simple. By definition there exists a point $y \in Y$ such that, if
$\{(H_i,\psi_i)\}_ {i=1,\ldots,s}$ is its combinatorial data, the
map $\nu: H=\oplus_i H_i \lrg G$ isn't injective. Therefore the
kernel $K$ of the associated combinatorial cover $\rho$ is not
trivial.

By hypothesis, $X$ is a locally complete intersection and so $X$
is Gorenstein over $y$.

\noindent Therefore by theorem \ref{theorema A in sezion proof},
the characters of the combinatorial data at the point $y$ of the
$\pi$-cover satisfy the identity (\ref{cond su caratteri in teoA
in sez proof}). Then the combinatorial data of the associated
combinatorial cover $\rho$ at the point $0 \in \mathbb{C}^s/H$ satisfy the
same identity. We studied the $(\ze_p)^n$-covers in example
\ref{esempio (Z_p)^n-cover con remark n_i mag 3} and remark
\ref{oss se gorens n_i div 1} shows that the codimension of the
singularities is at least three. Hence applying theorem
\ref{teorema schlessinger} as in remark \ref{osserv(Z_2)^3-cover
no loc sem no inters compl},  $\mathbb{C}^s/K$ is infinitesimally  rigid
and so it isn't a complete intersection. Applying theorem \ref{teo
inters completa se base e fibre}, the fiber $\rho^{-1}(0)$ is not
a locally complete intersection and so, by lemma \ref{lemm esiste
combina cover stesse fibre}, the fiber of $\pi$ over $y$ is not  a
locally complete intersection. Therefore, applying again theorem
\ref{teo inters completa se base e fibre}, $X$ is not a complete
intersection in $y$ and this is a contradiction.
\end{proof}

\begin{oss}
At this stage we could try to extend the previous equivalence to
each $G$-cover. As already observed, the implication i)
$\Rightarrow$ ii) still holds for each $G$-cover (remark \ref{oss
loc siple implica sempre loc compl inters}). In section \ref{sezio
key example}, example \ref{example Z_pqr cover} shows the
existence of a local complete intersection $G$-cover that isn't
locally simple and so ii) $\Rightarrow$ i) doesn't hold in
general. Probably for the general case, we have to generalize the
notion of locally simple cover to have an equivalence with locally
complete intersection.
\end{oss}

\bigskip

\noindent {\bf Acknowledgments.} I wish to thank M. Manetti for
having introduced me to this problem and for all remarks,
suggestions and encouragements; several ideas of this paper are
grown under his influence. I am very grateful to the referees for
their improvements in the presentation of the paper.

\end{document}